\documentclass[12 pt,a4paper,fleqn]{article}
\usepackage[latin1]{inputenc}
\usepackage{amsfonts}
\usepackage{amsmath}
\usepackage{amssymb}
\usepackage{graphicx}
\usepackage{hyperref}
\usepackage[left=2.5cm, right=2.5cm, bottom=3cm]{geometry}

\begin{document}

\title{About a new family of sequences}
	\author{Felipe Bottega Diniz}
	\maketitle
		
\begin{quote}
\begin{center}
	{\bf Abstract} 
\end{center}
First we define a new kind of function over $\mathbb{N}$. For each $i\in\mathbb{N}$ we
have an associated function, which will be called $S_i$ . Then we define a new kind of
sequence, to be made from the functions $S_i$ . Finally, we will see that some of these
sequences has a self-similarity feature.
\end{quote}	

\section{Introduction}
	There are sequences that are of great value for mathematics, many sequences are
also important for other sciences. For these reasons the study of sequences is very
important. But we can not forget that creativity in mathematics is also important,
and this is more true when we are dealing with sequences, it is not enough just to
study sequences, we also need to create sequences.
This work deals with a new kind of sequence, we will see that this kind of sequence
has interesting properties, that is reason enough to make this work worthwhile.

	For this work, we are considering $\mathbb{N} = \{1, 2, 3, \ldots\}$.
Consider $t \in \mathbb{N}$, such that $t$ has $n \in \mathbb{N}$ digits, so we can write $t = a_n a_{n-1} \ldots a_2 a_1$.
We can divide $t$ in blocks of $i$ digits (except at most one block), from right to left,
with $i \in \mathbb{N}$ and $i \leq n$. If $i$ divides $n$, then all blocks has $i$ digits, otherwise the last block will have
fewer digits.

	Below are some examples to get a better understanding. Consider the number $123456$,
we divide it into blocks of $i$ digits. Also consider that each block will be bracketed.
number of digits.\\\\
	${}\hspace{5cm}$Blocks of 1 digit: $[1] [2] [3] [4] [5] [6]$\\
	${}\hspace{5cm}$Blocks of 2 digits:$[12] [34] [56]$\\
	${}\hspace{5cm}$Blocks of 3 digits:$[123] [456]$\\
	${}\hspace{5cm}$Blocks of 4 digits:$[12] [3456]$\\
	${}\hspace{5cm}$Blocks of 5 digits:$[1] [23456]$\\
	${}\hspace{5cm}$Blocks of 6 digits:$[123456]$\\
	${}\hspace{5cm}$Blocks of 7 digits:$[123456]$\\
	${}\hspace{7cm}\vdots$
		
	Define the function $T:\mathbb{N}\to\mathbb{N}$ such that $T(m) = ``\text{Sum of digits of } m"\cdot m$. Consider the number $t$ divided into blocks of $i$ digits. Let $[a_k a_{k-1} \ldots a_{k-i+1}]$ be one of its blocks, then $T([a_k a_{k-1} \ldots a_{k-i+1}]) = (a_k + a_{k-1} +\ldots + a_{k-i+1})\cdot (a_k a_{k-1} \ldots a_{k-i+1})$. With this in mind, now we define the function $S_i:\mathbb{N}\to\mathbb{N}$ to be $$ S_i(t) = \sum T([a_k a_{k-1} \ldots a_{k-i+1}]), $$
where the sum goes over all the blocks as we described before.\\

	Taking our previous example $t = 123456$, we have the following.
	
	$$S_1 (123456) = T (1) + T (2) + T (3) + T (4) + T (5) + T (6) =
 1 \cdot 1 + 2 \cdot 2 + 3 \cdot 3 + 4 \cdot 4 + 5 \cdot 5 + 6 \cdot 6 $$ 
$$= 1 + 4 + 9 + 16 + 25 + 36 = 91$$

	$$S_2 (123456) = T (12) + T (34) + T (56) = (1 + 2) \cdot 12 + (3 + 4) \cdot 34 + (5 + 6) \cdot 56 
= 3 \cdot 12 + 7 \cdot 34 + 11 \cdot 56 = 890$$

	$$S_3 (123456) = T (123) + T (456) = (1 + 2 + 3) \cdot 123 + (4 + 5 + 6) \cdot 456 
= 6 \cdot 123 + 15 \cdot 456 = 7578$$

	$$S_4 (123456) = T (12) + T (3456) = (1 + 2) \cdot 12 + (3 + 4 + 5 + 6) \cdot 3456 
= 3 \cdot 12 + 18 \cdot 3456 = 62244$$

	$$S_5 (123456) = T (1) + T (23456) = 1 \cdot 1 + (2 + 3 + 4 + 5 + 6) \cdot 23456 
= 1 + 20 \cdot 23456 = 469121$$

	$$S_6 (123456) = T (123456) = (1 + 2 + 3 + 4 + 5 + 6) \cdot 123456 
= 21 \cdot 123456 = 2592576$$

	$$S_7 (123456) = T (123456) = (1 + 2 + 3 + 4 + 5 + 6) \cdot  123456 
= 21 \cdot 123456 = 2592576	$$

	$$\vdots $$
	
	We can also define $S_\infty(t)$, which of course means $``\text{Sum of digits of } t"\cdot t$. Note that this
definition of $S_\infty (t)$ coincides with the definition of the sequence A057147, hence, this
characterization represents a more general view.

\section{Some Results}
		
	\textbf{Theorem 1: } If $0\leq t\leq 9$, then $S_i(t) = t^2$, for all $i\in\mathbb{N}$.\\

	\textbf{Proof: } If $1 \leq i \leq 9$, then it's clear that $S_i(t) = (t) \cdot t = t^2.$\\\\	
	\textbf{Theorem 2: } Let $k, t \in \mathbb{N}$ such that $t = a_n a_{n-1} \ldots a_1$ and $k$ divides $a_m$ for each
$m \in \{1, 2, \ldots , n\}$. Then, for all $i \in \mathbb{N}$, exists $b \in \mathbb{N}$ such that $S_i (t) = k^2\cdot S_i (b)$.\\

	\textbf{Proof: } Since $k$ divides all the digits of $t$, for each digit $a_m$ of $t$ there is a natural
$b_m$ such that $a_m = k \cdot b_m$. Thus, we have that $t = a_n a_{n-1} \ldots a_1 = k \cdot b_n\ k \cdot b_{n-1}\ k \cdot b_1$. 
$$S_i(t) = S_i(k \cdot b_n\ k \cdot b_{n-1}\  k \cdot b_1 ) = $$
$$= S_i(k\cdot b_n\ k\cdot b_{n-1}\ \ldots k\cdot b_{n-q} )+ \ldots + S_i(k\cdot b_i\ \ldots k\cdot b_1 ) = $$
$$ =(k\cdot b_n +k\cdot b_{n-1} + \ldots +k\cdot b_{n-q} )\cdot (k\cdot b_n\ \ k\cdot b_{n-1} \ldots k\cdot b_{n-q} )+ \ldots +(k\cdot b_i + \ldots +k\cdot b_1 )\cdot(k\cdot b_i \ldots k\cdot b_1 ) = $$
$$= k(b_n +b_{n-1} + \ldots +b_{n-q} )\cdot(k\cdot b_n\ k\cdot b_{n-1} \ldots k\cdot b_{n-q} )+ \ldots +k(b_i + \ldots +b_1 )\cdot (k\cdot b_i \ldots k\cdot b_1 ) =$$
$$= k \Big((b_n +b_{n-1} + \ldots +b_{n-q} )\cdot(k\cdot b_n\ k\cdot b_{n-1} \ldots k\cdot b_{n-q} )+ \ldots +(b_i + \ldots +_b1 )\cdot (k\cdot b_i \ldots k\cdot b_1 ) \Big).$$

	Note that any number of the form $k \cdot b_x\ k \cdot b_{x-1}  \ldots k \cdot b_{x-y}$ can be written as
$$10^y \cdot k \cdot b_x + 10^{y-1} \cdot k \cdot b_{x-1} + \ldots + 10^0 \cdot k \cdot b_{x-y}.$$ We can write the last equation as
follows:

$$\hspace{-1cm} k \Big((b_n +b_{n-1} + \ldots +b{n-q} )\cdot (10^q \cdot k\cdot b_n +10^{q-1} \cdot k\cdot b_{n-1} + \ldots +10^0 \cdot k\cdot b_{n-q} )+ \ldots +(b_i + \ldots +b_1 )\cdot (10^i \cdot k\cdot b_i \ldots 10^0 \cdot k\cdot b_1 )\Big) = $$
$$ = k \Big((b_n +b_{n-1} + \ldots +b_{n-q} )\cdot k(10^q \cdot b_n +10^{q-1} \cdot b_{n-1} + . . . +10^0 \cdot b_{n-q} )+ \ldots +(b_i + \ldots+b_1 )\cdot k(10^i \cdot b_i \ldots 10^0 \cdot b_1 )\Big) = $$
$$= k^2 \Big((b_n +b_{n-1} + \ldots +b_{n-q} )\cdot(10^q \cdot b_n +10^{q-1} \cdot b_{n-1} + \ldots +10^0 \cdot b_{n-q} )+ \ldots +(b_i + \ldots +b_1 )\cdot (10^i \cdot b_i \ldots 10^0 \cdot b_1 )\Big) =$$
$$= k^2 \Big((b_n +b_{n-1} + \ldots +b_{n-q} )\cdot (b_n b_{n-1} \ldots b_{n-q} )+ \ldots +(b_i + \ldots + b_1 )\cdot (b_i \ldots b_1 )\Big).$$

	Let $b \in \mathbb{N}$ such that $b = b_n b_{n-1} \ldots b_1$, the existence of $b$ is guaranteed, because it is
guaranteed the existence of each of $b_m$ . Also note that the above equation is equal to $k^2 \cdot S_i(b)$, as we wanted to prove.\\\\
	\textbf{Corollary 2.1: } If $0\leq a \leq 9$, then for each $i \in \mathbb{N}$, $S_i (\underbrace{a\ a \ldots a}_{n\ times }) = a^2 \cdot S_i (\underbrace{1\ 1 \ldots 1}_{n\ times})$.\\
	
	\textbf{Proof: } Note that $(a\ a \ldots a) = (a\cdot1\ a\cdot 1 \ldots a\cdot 1)$. Using theorem 2, we get $S_i(a\cdot1\ a\cdot 1 \ldots a\cdot 1) = a^2\cdot S_i(1\ 1\ \ldots 1)$, as desired.\\\\
	\textbf{Corollary 2.2: } If $S_i(t)$ is a prime number, then there is not a number (greater
than 1) that divides all the digits of t.\\

	\textbf{Proof: } By Theorem 2, if all the digits of t are divisible by a number (greater than 1),
then $S_i(t)$ can't be prime because it's the product of two numbers greater than 1.\\\\
	\textbf{Theorem 3: } No function $S_i$ is injective.\\

	\textbf{Proof: } For every $i\in\mathbb{N}$ , we have that $S_i (1) = 1$, we also have that $S_i(10^i) = S_i(\underbrace{100\ldots 0}_{i+1\ digits}) = 1$. Therefore, $S_i$ is not injective.\\\\
	\textbf{Theorem 4: } Every function $S_i$ is surjective.\\

	\textbf{Proof: } Let $n, t \in \mathbb{N}$ such that $n$ is arbitrary and $t=1\underbrace{\underbrace{0\ 0\ \ldots 1}_{i\ digits}\ldots\underbrace{0\ 0\ \ldots 1}_{i\ digits}}_{n-1\ times}$. Then $S_i(t) = 1 + \underbrace{1+1+\ldots +1}_{n-1\ times} = n$. Therefore, $S_i$ is surjective.\\\\
	
	We can construct some special sequences from the functions $S_i$. For each $i \in \mathbb{N}$, we will be interested in the sequence $\mathcal{S}_i = (S_i(1), S_i(2), S_i(3), \ldots)$.
	
\section{Charts}

	Below are some charts of some of these sequences. These charts clearly show that there is a kind of pattern. In particular, the charts of $S_1$ , $S_2$ and $S_7$ has some self-similarity features.\footnotesize
	
	\begin{picture}(1,1)
		\put(60,-220){\includegraphics[scale=0.8]{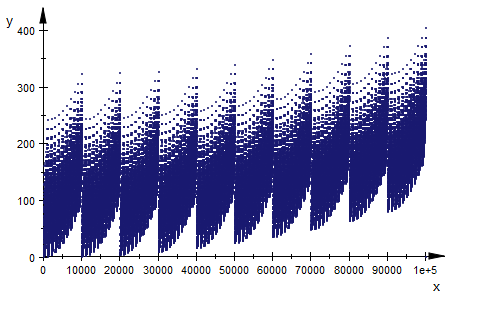}}
	\end{picture}\\\\\\\\\\\\\\\\\\\\\\\\\\\\\\\\\\\\
	\begin{center}
		Fig. 1: values of $\mathcal{S}_1(x)$, with $x=1\ldots 10^5$. 	
	\end{center}\newpage	
	
	\begin{picture}(1,1)
		\put(60,-220){\includegraphics[scale=0.8]{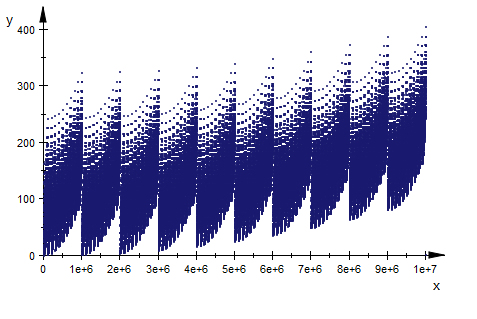}}
	\end{picture}\\\\\\\\\\\\\\\\\\\\\\\\\\\\\\\\\\\\
	\begin{center}
		Fig. 2: values of $\mathcal{S}_1(x)$, with $x=1\ldots 10^7$. 	
	\end{center}
	
	\begin{picture}(1,1)
		\put(60,-220){\includegraphics[scale=0.8]{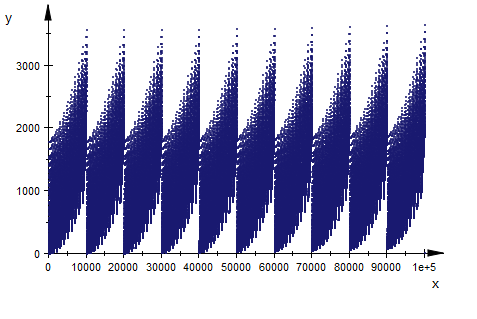}}
	\end{picture}\\\\\\\\\\\\\\\\\\\\\\\\\\\\\\\\\\\\
	\begin{center}
		Fig. 3: values of $\mathcal{S}_2(x)$, with $x=1\ldots 10^5$. 	
	\end{center}\newpage
	
	\begin{picture}(1,1)
		\put(60,-220){\includegraphics[scale=0.8]{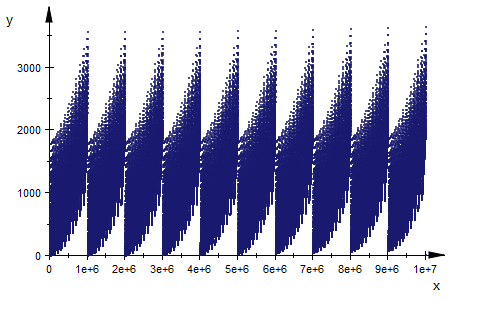}}
	\end{picture}\\\\\\\\\\\\\\\\\\\\\\\\\\\\\\\\\\\\
	\begin{center}
		Fig. 4: values of $\mathcal{S}_2(x)$, with $x=1\ldots 10^7$. 	
	\end{center}
	
	\begin{picture}(1,1)
		\put(60,-220){\includegraphics[scale=0.8]{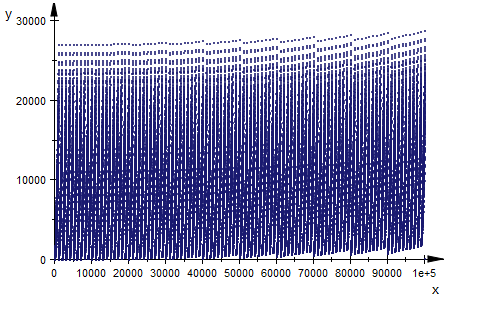}}
	\end{picture}\\\\\\\\\\\\\\\\\\\\\\\\\\\\\\\\\\\\
	\begin{center}
		Fig. 5: values of $\mathcal{S}_3(x)$, with $x=1\ldots 10^5$. 	
	\end{center}\newpage
	
	\begin{picture}(1,1)
		\put(60,-220){\includegraphics[scale=0.8]{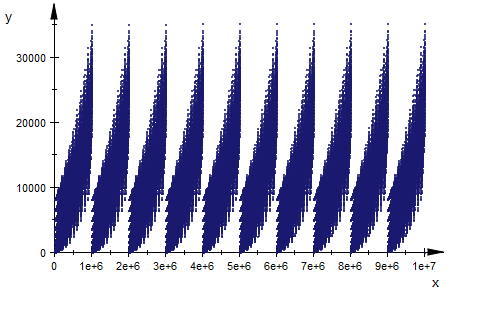}}
	\end{picture}\\\\\\\\\\\\\\\\\\\\\\\\\\\\\\\\\\\\
	\begin{center}
		Fig. 6: values of $\mathcal{S}_3(x)$, with $x=1\ldots 10^7$. 	
	\end{center}
	
	\begin{picture}(1,1)
		\put(60,-220){\includegraphics[scale=0.8]{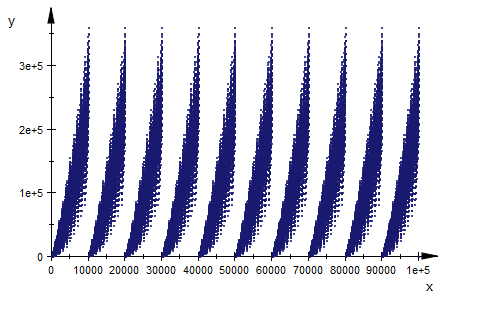}}
	\end{picture}\\\\\\\\\\\\\\\\\\\\\\\\\\\\\\\\\\\\
	\begin{center}
		Fig. 7: values of $\mathcal{S}_4(x)$, with $x=1\ldots 10^5$. 	
	\end{center}\newpage
	
	\begin{picture}(1,1)
		\put(60,-220){\includegraphics[scale=0.8]{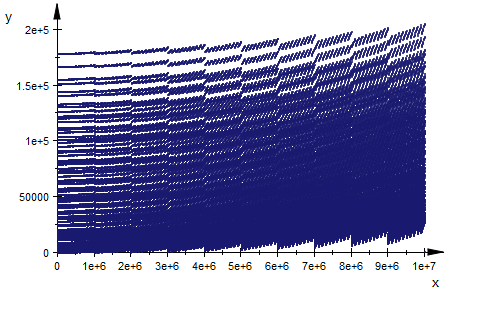}}
	\end{picture}\\\\\\\\\\\\\\\\\\\\\\\\\\\\\\\\\\\\
	\begin{center}
		Fig. 8: values of $\mathcal{S}_4(x)$, with $x=1\ldots 10^7$. 	
	\end{center}
	
	\begin{picture}(1,1)
		\put(60,-220){\includegraphics[scale=0.8]{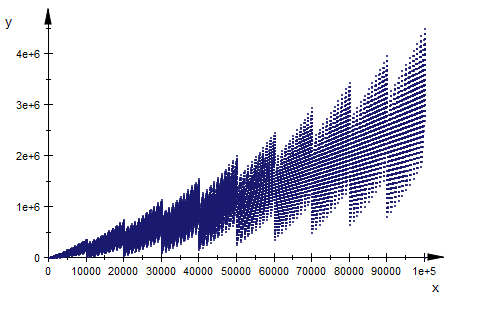}}
	\end{picture}\\\\\\\\\\\\\\\\\\\\\\\\\\\\\\\\\\\\
	\begin{center}
		Fig. 9: values of $\mathcal{S}_5(x)$, with $x=1\ldots 10^5$. 	
	\end{center}\newpage
	
	\begin{picture}(1,1)
		\put(60,-220){\includegraphics[scale=0.8]{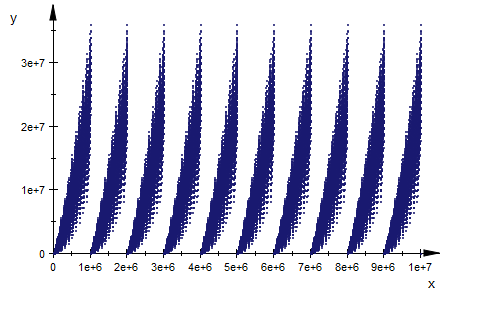}}
	\end{picture}\\\\\\\\\\\\\\\\\\\\\\\\\\\\\\\\\\\\
	\begin{center}
		Fig. 10: values of $\mathcal{S}_6(x)$, with $x=1\ldots 10^7$. 	
	\end{center}
	
	\begin{picture}(1,1)
		\put(60,-220){\includegraphics[scale=0.8]{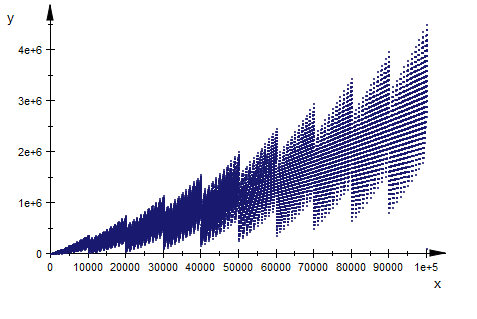}}
	\end{picture}\\\\\\\\\\\\\\\\\\\\\\\\\\\\\\\\\\\\
	\begin{center}
		Fig. 11: values of $\mathcal{S}_7(x)$, with $x=1\ldots 10^5$. 	
	\end{center}\newpage
	
	\begin{picture}(1,1)
		\put(60,-220){\includegraphics[scale=0.8]{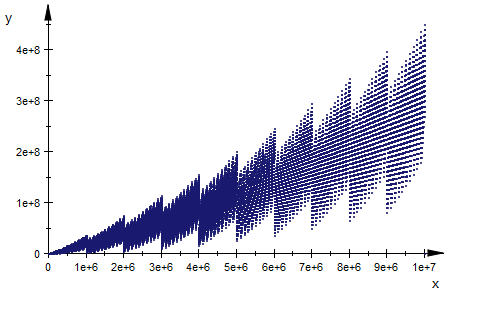}}
	\end{picture}\\\\\\\\\\\\\\\\\\\\\\\\\\\\\\\\\\\\
	\begin{center}
		Fig. 12: values of $\mathcal{S}_7(x)$, with $x=1\ldots 10^7$. 	
	\end{center}\newpage
	
\section*{References}	

	[1] Fabio B. Martinez, Carlos G. Moreira, Nicolau Saldanha and Eduardo Tengan. Teoria dos Numeros, IMPA, 2010.\\\\
	{}[2] Elon L. Lima. Curso de Analise Volume 1, IMPA, 2009.\\ \\
	{}[3] Jose P. Santos. Introducao a Teoria dos Numeros, IMPA, 2009.\\\\ 
	{}[4] \href{https://oeis.org/A057147}{https://oeis.org/A057147}.\\\\
	{}[5] \href{https://oeis.org/A213630}{https://oeis.org/A213630}.

\end{document}